\documentclass[a4paper,12pt]{article}
\usepackage{amsmath}
\usepackage{color}
\usepackage{amssymb}
\usepackage{latexsym}
\usepackage{url}
\numberwithin{equation}{section}

\def\R{{\bf R}}

\def\N{{\bf N}}

\def\d{\displaystyle}
\def\e{{\varepsilon}}

\def\wt{\widetilde}

\newtheorem{thm}{Theorem}[section]

\newtheorem{lem}{Lemma}[section]

\newtheorem{rem}{Remark}[section]
\newtheorem{Def}{Definition}[section]

\title{Nonexistence of global solutions
of nonlinear wave equations with weak time-dependent damping
related to Glassey conjecture
}
\author{
Ning-An Lai
\footnote{Department of Mathematics,
Lishui University,
Lishui City 323000,
China.
e-mail: hyayue@gmail.com.}
\quad
Hiroyuki Takamura
\footnote{Department of Complex and Intelligent Systems,
Faculty of Systems Information Science,
Future University Hakodate,
116-2 Kamedanakano-cho,
Hakodate, Hokkaido 041-8655, Japan.
e-mail: takamura@fun.ac.jp.}
}
\date{
\[
\begin{array}{ll}
\mbox{\footnotesize{\bf Keywords:}}
& \mbox{\footnotesize weak damping, nonlinear wave equations, blow-up, lifespan}\\
\mbox{\footnotesize{\bf MSC2010:}}
& \mbox{\footnotesize primary 35L71, secondary 35B44}\\
\end{array}
\]
}

\begin{document}
\maketitle
\begin{abstract}
This work is devoted to the nonexistence of global-in-time
energy solutions of nonlinear wave equation of derivative type
with weak time-dependent damping in the scattering and scale invariant range.
By introducing some multipliers to absorb the damping term,
we succeed in establishing the same upper bound of the lifespan for the scattering damping
as the non-damped case,
which is a part of so-called Glassey conjecture on nonlinear wave equations.
We also study an upper bound of the lifespan for the scale invariant damping
with the same method.
\end{abstract}


\section{Introduction}
\par\quad
In this work, we consider the following Cauchy problem for
the nonlinear damped wave equations.
\begin{equation}
\label{IVP}
\left\{
\begin{array}{l}
\d u_{tt}-\Delta u+\frac{\mu}{(1+t)^\beta}u_t=|u_t|^p
\quad \mbox{in}\ \R^n\times[0,\infty),\\
u(x,0)=\e f(x),\ u_t(x,0)=\e g(x), \quad  x\in\R^n,
\end{array}
\right.
\end{equation}
where $\mu\ge0,\ n\in\N$ and $\beta\geq 1$.
We assume that $\e>0$ is a \lq\lq small" parameter
and that $f,g$ are in the energy space with compact support.
The restriction on $\beta$ is so-called scattering case ($\beta>1$) in which the solution of the linear equation
scatters to the one of free wave equations, and scale invariant case ($\beta=1$), in which the linear equation in \eqref{IVP}
is invariant under the following scaling transform
\[
\wt{u}(x,t):=u(\sigma x, \sigma(1+t)-1),\ \sigma>0.
\]
We refer the reader to Wirth \cite{Wirth1, Wirth2, Wirth3}
for the classifications on $\beta$.
\par
First we shall outline the results on (\ref{IVP}) without damping, i.e. $\mu=0$.
It has been conjectured that there is no global solution
for $p>1$ when $n=1$,
and also that there is a critical power
\[
p_c(n):=\frac{n+1}{n-1}
\]
in the sense that we have global existence for $p>p_c(n)$
while the blow-up in finite time occurs for $1<p\le p_c(n)$
when $n\ge2$.
This problem is so-called Glassey conjecture appeared in Glassey \cite{Glassey},
and was initiated by John \cite{John1}
in which he studied more general equations for $n=3$,
and proved that the solution blows-up for $p=2$.
We note that his method works also for $1<p\le2$.
After \cite{John1},
Masuda \cite{Masuda} obtained the blow-up result for $p=2$ and $n=1, 2, 3$.
Schaeffer \cite{Schaeffer} established a blow-up result
for $n=2$ and $p=3$,
and conjectured that $p_c(2)=3$.
See also John \cite{John2}.
Agemi \cite{Agemi} extended the result in \cite{Schaeffer} to $1<p\le3$.
Moreover, Rammaha \cite{Rammaha} studied the blow-up result
for high dimensional case, $n\ge4$, under the radially symmetric assumption.
Finally, Zhou \cite{Zhou1} introduced a simple proof of the blow-up result for
all $n\ge2$ and $1<p\leq p_c(n)$
as well as $p>1$ and $n=1$,
and obtained the upper bound of lifespan
of the solution.
For global existence part, Sideris \cite{Sideris} proved it for $n=3$ and $p\geq 2$
under the radially symmetric assumption.
Hidano and Tsutaya \cite{Hidano1},
and independently Tzvetkov \cite{Tzvetkov},
obtained the global-in-time solution for $n=2, 3$ and $p>p_c(n)$
without radially symmetric assumption.
Finally, Hidano, Wang and Yokoyama \cite{Hidano2} generalized
the global existence result to high dimensional cases, i.e. $n\ge4$,
under the radially symmetric assumption.
\par
In this work, we are going to study Cauchy problem (\ref{IVP}) for $\mu>0$.
We focus on the blow-up result and lifespan estimate from above.
Without the damping term, the corresponding results has been obtained
in Zhou \cite{Zhou1}, as mentioned above.
For our problem we have to overcome the difficulty caused by the damping term. However, due to the scattering and scale invariant coefficients,
we may use the multipliers introduced in the authors \cite{LT} and the authors and Wakasa \cite{LTW} respectively,
to absorb the damping term.
Then by combining the method used in Zhou and Han \cite{Zhou2},
we get the blow-up result and the upper bound of lifespan estimate.

\section{Main Result}
\par\quad
Before showing the main result, we first define the energy and weak solution of the Cauchy problem (\ref{IVP}).
\begin{Def}\label{def1}
 As in \cite{LT} and \cite{LTW}, we say that $u$ is an energy solution of
 (\ref{IVP}) on $[0,T)$
if
\[
u\in C([0,T),H^1(\R^n))\cap C^1([0,T),L^2(\R^n))
\cap C^1((0,T),L^p(\R^n)),
\]
satisfies
\begin{equation}
\label{energysol}
\begin{array}{l}
\d\quad\int_{\R^n}u_t(x,t)\phi(x,t)dx-\int_{\R^n}u_t(x,0)\phi(x,0)dx\\
\d\quad+\int_0^tds\int_{\R^n}\left\{-u_t(x,s)\phi_t(x,s)+\nabla u(x,s)\cdot\nabla\phi(x,s)\right\}dx\\
\d\quad+\int_0^tds\int_{\R^n}\frac{\mu u_t(x,s)}{(1+s)^{\beta}}\phi(x,s)dx\\
\d=\int_0^tds\int_{\R^n}|u_t(x,s)|^p\phi(x,s)dx
\end{array}
\end{equation}
with any $\phi\in C_0^{\infty}(\R^n\times[0,T))$ and any $t\in[0,T)$.
\end{Def}

Employing the integration by parts in \eqref{energysol}
and letting $t\rightarrow T$, we have that
\[
\begin{array}{l}
\d\quad\int_{\R^n\times[0,T)}
u(x,s)\left\{\phi_{tt}(x,s)-\Delta \phi(x,s)
-\left(\frac{\mu\phi(x,s)}{(1+s)^{\beta}}\right)_s \right\}dxds\\
\d=\int_{\R^n}\mu u(x,0)\phi(x,0)dx-\int_{\R^n}u(x,0)\phi_t(x,0)dx\\
\d\quad+\int_{\R^n}u_t(x,0)\phi(x,0)dx+\int_{\R^n\times[0,T)}|u_t(x,s)|^p\phi(x,s)dxds.
\end{array}
\]
This is exactly the definition of the weak solution of \eqref{IVP}.
\par
Our main results are stated in the following tow theorems.
\begin{thm}
\label{blowup}
Let $\mu>0$ and $\beta>1$.
Assume that both $f\in H^1(\R^n)$ and $g\in L^2(\R^n)$ are non-negative,
and $g$ does not vanish identically.
Suppose that an energy solution $u$ of \eqref{IVP} on $[0,T)$ satisfies
\begin{equation}
\label{support}
\mbox{\rm supp}\ u\ \subset\{(x,t)\in\R^n\times[0,\infty)\ :\ |x|\le t+R\}
\end{equation}
with some $R\ge1$.
Then, there exists a constant $\e_0=\e_0(f,g,n,p,\mu, \beta, R)>0$
such that $T$ has to satisfy
\[
T\leq
\left\{
\begin{array}{ll}
 C\e^{-(p-1)/\{1-(n-1)(p-1)/2\}}
 &
 \mbox{for}\
 1<p<
 \left\{
 \begin{array}{ll}
 p_c(n) & \mbox{when}\ n\ge2,\\
 \infty & \mbox{when}\ n=1,
 \end{array}
 \right.
 \\
 \exp\left(C\e^{-(p-1)}\right)
&
\mbox{for}\ p=p_c(n)\ \mbox{and}\ n\geq 2
\end{array}
\right.
\]
with $0<\e\le\e_0$, where $C$ is a positive constant independent of $\e$.
\end{thm}

\begin{rem}
This estimate provides us the same upper bound of the lifespan
as the case of $\mu=0$ in Zhou \cite{Zhou1}.
\end{rem}

\begin{thm}
\label{blowup1}
Let $\mu>0$ and $\beta=1$.
Assume the same condition on $f,g$ and supp $u$
to Theorem \ref{blowup}.
Then, for $n\geq 1$, there exists a constant\\
$\e_0=\e_0(f,g,n,p,\mu, \beta, R)>0$
such that $T$ has to satisfy
\[
T\leq
\left\{
\begin{array}{ll}
 C\e^{-(p-1)/\{1-(n+2\mu-1)(p-1)/2\}}
 &
 \mbox{for}\
 1<p<p_c(n+2\mu),
 \\
 \exp\left(C\e^{-(p-1)}\right)
&
\mbox{for}\ p=p_c(n+2\mu)
\end{array}
\right.
\]
with $0<\e\le\e_0$, where $C$ is a positive constant independent of $\e$.
\end{thm}

\begin{rem}
Along with the definition of the scattering case by Wirth \cite{Wirth1, Wirth2, Wirth3},
Theorem \ref{blowup} can be established for generalized damping
for which $\mu/(1+t)^\beta u_t$ in (\ref{IVP}) is replaced
by positive function $b(t)u_t$ satisfying $b\in L^1([0,\infty))$.
It is easy to prove this fact by our proof below if one substitutes
the definition of the multiplier $m$ in (\ref{test1}) by
\[
m(t)=\exp\left(-\int_t^\infty b(s)ds\right)
\]
due to the fact that we only need a boundedness of $m$.
But such a generalization can not be available in Theorem 2.2
due to the unboundedness of the multiplier $m_1(t)$ in (\ref{test1in}) below.
\end{rem}

\section{Proof of Theorem \ref{blowup}}
\par\quad
In the proof of the main theorem,
we make use of two key tools.
The first one is a multiplier
\begin{equation}
\label{test1}
\begin{aligned}
m(t):=\exp\left(\mu\frac{(1+t)^{1-\beta}}{1-\beta}\right),
\end{aligned}
\end{equation}
which was first introduced in Lai and Takamura \cite{LT} and has a property
\[
\frac{m'(t)}{m(t)}=\frac{\mu}{(1+t)^\beta}.
\]
This multiplier is specially useful for the study
of nonlinear damped wave equation
with $\beta>1$ due to its boundedness
from above and below as
\begin{equation}
\label{bound}
1\ge m(t) \ge m(0)>0 \quad\mbox{for}\ t\ge0.
\end{equation}
The other one is defined as
\begin{equation}
\label{test11}
\psi(x,t):=e^{-t}\phi_1(x),
\quad
\phi_1(x):=
\left\{
\begin{array}{ll}
\d\int_{S^{n-1}}e^{x\cdot\omega}dS_\omega & \mbox{for}\ n\ge2,\\
e^x+e^{-x} & \mbox{for}\ n=1,
\end{array}
\right.
\end{equation}
which was introduced in Yordanov and Zhang \cite{Yordanov}
and admits the following good properties:
\begin{equation}
\label{psi}
\psi_t=-\psi,\quad\psi_{tt}=\Delta\psi=\psi.
\end{equation}
We note that there exists a constant $C_1=C_1(n,R)>0$ such that
\begin{equation}
\label{inequalityforpsi}
\int_{|x|\leq t+R}\psi(x, t)dx \leq C_1(t+1)^{(n-1)/2}
\quad\mbox{for}\ t\ge0
\end{equation}
with any constant $R>0$.
\par
Setting
\begin{equation}
\label{F1def}
F_1(t):=\int_{\R^n}u(x, t)\psi(x, t)dx,
\end{equation}
we have the following lemma.
\begin{lem}
\label{F1}
Under the same assumption of Theorem \ref{blowup}, it holds that
\begin{equation}
\label{F1postive}
F_1(t)\ge \frac{m(0)\e}{2}\int_{\R^n}f(x)\phi_1(x)dx\ge0
\quad\mbox{for}\ t\ge0.
\end{equation}
\end{lem}
{\it Proof.} The proof of Lemma \ref{F1} is parallel to that of (3.9) in \cite{LT}. For convenience, we write down the details. By the definition \eqref{energysol}, we get
\[
\begin{array}{l}
\d\frac{d}{dt}\int_{\R^n}u_t(x,t)\phi(x,t)dx
+\int_{\R^n}\frac{\mu u_t(x,t)}{(1+t)^\beta}\phi(x,t)dx\\
\d\quad+\int_{\R^n}\left\{-u_t(x,t)\phi_t(x,t)-u(x,t)\Delta\phi(x,t)\right\}dx\\
\d=\int_{\R^n}|u_t(x,t)|^p\phi(x,t)dx.
\end{array}
\]
Multiplying the both sides of the above equality by $m(t)$ we have
\[
\begin{array}{l}
\d\frac{d}{dt}\left\{m(t)\int_{\R^n}u_t(x,t)\phi(x,t)dx\right\}\\
\d\quad+m(t)\int_{\R^n}\left\{-u_t(x,t)\phi_t(x,t)-u(x,t)\Delta\phi(x,t)\right\}dx\\
\d=m(t)\int_{\R^n}|u_t(x,t)|^p\phi(x,t)dx.
\end{array}
\]
Integration this equality over $[0,t]$ implies that
\[
\begin{array}{l}
\d m(t)\int_{\R^n}u_t(x,t)\phi(x,t)dx
-m(0)\e\int_{\R^n}g(x)\phi(x,0)dx\\
\d\quad-\int_0^tm(s)ds\int_{\R^n}\left\{
u_t(x,s)\phi_t(x,s)+u(x,s)\Delta\phi(x,s)\right\}dx\\
\d=\int_0^tm(s)ds\int_{\R^n}|u_t(x,s)|^p\phi(x,s)dx.
\end{array}
\]
Replacing $\phi(x, t)$ with $\psi(x, t)$ on supp $u$ in the above inequality,
making use of (\ref{psi}) and integration by parts in $t$-integral
in the second line, we come to
\[
\begin{array}{l}
\d m(t)\{F_1'(t)+2F_1(t)\}-m(0)\e\int_{\R^n}\left\{f(x)+g(x)\right\}\phi_1(x)dx\\
\d=\int_0^tm'(s)F_1(s)ds+\int_0^tm(s)ds\int_{\R^n}|u_t(x,s)|^p\psi(x,s)dx.
\end{array}
\]
which yields
\begin{equation}
\label{F1inequality}
\begin{aligned}
F'_1(t)+2F_1(t)
&\d \ge\frac{m(0)}{m(t)}
C_{f,g}\e+\frac{1}{m(t)}\int_0^tm(s)
\frac{\mu}{(1+s)^\beta}F_1(s)ds\\
&\d \ge m(0)C_{f,g}\e+\frac{1}{m(t)}\int_0^tm(s)
\frac{\mu}{(1+s)^\beta}F_1(s)ds,
\end{aligned}
\end{equation}
where
\[
C_{f,g}:=\int_{\R^n}\left\{f(x)+g(x)\right\}\phi_1(x)dx.
\]
Here we have used the boundedness of $m$ in (\ref{bound}).
Hence it is easy to get from \eqref{F1inequality} that
\[
\begin{array}{ll}
e^{2t}F_1(t)
&\d\ge F_1(0)+m(0)C_{f,g}\e\int_0^te^{2s}ds\\
&\d\quad+\int_0^t\frac{e^{2s}}{m(s)}ds
\int_0^sm(r)\frac{\mu}{(1+r)^\beta}F_1(r)dr,\\
\end{array}
\]
which leads, by comparison argument, to
\[
e^{2t}F_1(t)\ge m(0)C_{f,0}\e+\frac{m(0)C_{f,0}\e}{2}(e^{2t}-1),
\]
and finally to
\begin{equation}
\label{18}
F_1(t)\ge\frac{1}{2}m(0)C_{f,0}\e
\quad\mbox{for}\ t\ge0
\end{equation}
which is exactly the inequality \eqref{F1postive} we need.
\hfill $\Box$
\vskip10pt
\par
Now we are in a position to prove Theorem \ref{blowup}. First we have
\begin{equation}
\label{u_t+u}
\begin{array}{l}
\d\frac{d}{dt}\left[m(t)\int_{\R^n}\left\{u_t(x, t)+u(x, t)\right\}\psi(x, t)dx\right]\\
\d=\frac{\mu}{(1+t)^\beta}m(t)\int_{\R^n}\left\{u_t(x, t)+u(x, t)\right\}\psi(x, t)dx\\
\quad\d+m(t)\frac{d}{dt}\int_{\R^n}\left\{u_t(x, t)+u(x, t)\right\}\psi(x, t)dx.
\end{array}
\end{equation}
Replacing the test function $\phi$ in the definition \eqref{energysol} by $\psi$
and taking derivative to both sides with respect to $t$, we have that
\begin{equation}
\label{derivateofenergysol}
\begin{array}{l}
\d\frac{d}{dt}\int_{\R^n}u_t(x,t)\psi(x,t)dx-\int_{\R^n}u_t(x,t)\psi_t(x,t)dx\\
\d+\int_{\R^n}\nabla u(x,t)\cdot \nabla\psi(x,t)dx
+\frac{\mu}{(1+t)^\beta}\int_{\R^n}u_t(x,t)\psi(x,t)dx\\
\d=\int_{\R^n}|u_t(x,t)|^p\psi(x,t)dx.
\end{array}
\end{equation}
Hence the integration by parts in the first term in the second line
of (\ref{derivateofenergysol}) with (\ref{psi}) yields that
\[
\begin{array}{l}
\d\frac{d}{dt}\int_{\R^n}\left\{u_t(x,t)+u(x,t)\right\}\psi(x,t)dx\\
\d=\int_{\R^n}|u_t(x,t)|^p\psi(x,t)dx
-\frac{\mu}{(1+t)^\beta}\int_{\R^n}u_t(x,t)\psi(x,t)dx.
\end{array}
\]
Plugging this equality into \eqref{u_t+u} we have
\begin{equation}
\label{derivateu_t+u}
\begin{array}{l}
\d\frac{d}{dt}\left[m(t)\int_{\R^n}\left\{u_t(x, t)+u(x, t)\right\}\psi(x, t)dx\right]\\
\d=m(t)\int_{\R^n}|u_t(x,t)|^p\psi (x,t)dx
+\frac{\mu}{(1+t)^\beta}m(t)F_1(t)
\end{array}
\end{equation}
for $t\ge0$.
Then (\ref{derivateu_t+u}) and the positivity of $F_1$ by Lemma \ref{F1} yield
\begin{equation}
\label{integrationu_t+u}
\begin{array}{l}
\d m(t)\int_{\R^n}\left\{u_t(x, t)+u(x, t)\right\}\psi(x, t)dx\\
\d\geq m(0)\e\int_{\R^n}\{f(x)+g(x)\}\phi_1(x)dx\\
\d\quad
+\int_0^tds\int_{\R^n}m(s)|u_t(x, s)|^p\psi(x, s) dx.
\end{array}
\end{equation}
\par
On the other hand,
\eqref{derivateofenergysol} also yields that
\[
\begin{array}{l}
\d\frac{d}{dt}\int_{\R^n}u_t(x,t)\psi(x,t)dx
+\frac{m'(t)}{m(t)}\int_{\R^n}u_t(x,t)\psi(x,t)dx\\
\d+\int_{\R^n}\left\{u_t(x,t)-u(x,t)\right\}\psi(x,t)dx\\
\d=\int_{\R^n}|u_t(x,t)|^p\psi(x,t)dx.
\end{array}
\]
Multiplying this equality by $m(t)$, we have
\begin{equation}
\label{derivateandintegration}
\begin{array}{l}
\d\frac{d}{dt}\left[m(t)\int_{\R^n}u_t(x,t)\psi(x,t)dx\right]\\
\d+m(t)\int_{\R^n}\left\{u_t(x,t)-u(x,t)\right\}\psi(x,t)dx\\
\d=m(t)\int_{\R^n}|u_t(x,t)|^p\psi(x,t)dx.
\end{array}
\end{equation}
Adding \eqref{integrationu_t+u} and \eqref{derivateandintegration} together,
we obtain that
\begin{equation}
\label{addtwo}
\begin{array}{l}
\d\frac{d}{dt}\left[m(t)\int_{\R^n}u_t(x,t)\psi(x,t)dx\right]
+2m(t)\int_{\R^n}u_t(x,t)\psi(x,t)dx\\
\d\geq m(0)\e\int_{\R^n}\{f(x)+g(x)\}\phi_1(x)dx
+m(t)\int_{\R^n}|u_t(x,t)|^p\psi(x,t)dx\\
\d\quad+\int_0^tm(s)ds\int_{\R^n}|u_t(x,s)|^p\psi(x,s)dx.\\
\end{array}
\end{equation}
Setting
\begin{equation}
\label{G}
\begin{array}{ll}
G(t):=
&
\d m(t)\int_{\R^n}u_t(x,t)\psi(x,t)dx
-\frac{m(0)\e}{2}\int_{\R^n}g(x)\phi_1(x)dx
\\
&
\d -\frac{1}{2}\int_0^tm(s)ds\int_{\R^n}|u_t(x,s)|^p\psi(x,s)dx,
\end{array}
\end{equation}
we have
\[
G(0)=\frac{m(0)\e}{2}\int_{\R^n}g(x)\phi_1(x)dx>0.
\]
Then it follows from (\ref{addtwo}) and direct computation that
\[
G'(t)+2G(t)
\geq \frac{m(t)}{2}\int_{\R^n}|u_t(x,t)|^p\psi(x,t)dx
+m(0)\e\int_{\R^n}\phi_1(x)f(x)dx
\geq 0
\]
which implies
\[
G(t)\geq e^{-2t}G(0)>0 \quad\mbox{for}\ t\ge0.
\]
Hence, by the definition \eqref{G}, it holds that
\begin{equation}
\label{inequalityforG(t)}
\begin{array}{l}
\d m(t)\int_{\R^n}u_t(x,t)\psi(x,t)dx\\
\d \ge\frac{1}{2}\int_0^tm(s)ds\int_{\R^n}|u_t(x,s)|^p\psi(x,s)dx\\
\d\quad+\frac{m(0)\e}{2}\int_{\R^n}g(x)\phi_1(x)dx.
\end{array}
\end{equation}
\par
Denoting
\[
H(t):=\frac{1}{2}\int_0^tm(s)ds\int_{\R^n}|u_t(x,s)|^p\psi(x,s)dx
+\frac{m(0)\e}{2}\int_{\R^n}g(x)\phi_1(x)dx,
\]
then, by \eqref{inequalityforG(t)}, we have
\begin{equation}
\label{inequalityformpsiu}
m(t)\int_{\R^n}u_t(x,t)\psi(x,t)dx\geq H(t) \quad\mbox{for}\ t\ge0.
\end{equation}
On the other hand, H\"{o}lder inequality and \eqref{inequalityforpsi} yield that
\begin{equation}
\label{inequalityfornonlinear}
2H'(t)\geq C_1^{1-p} (t+1)^{-(n-1)(p-1)/2}
\left(m(t)\int_{\R^n}u_t(x,t)\psi(x,t)dx\right)^p.
\end{equation}
Here we have used the boundedness of $m$ in (\ref{bound}) as $m(t)^{1-p}\ge1$.
We then conclude by \eqref{inequalityformpsiu} and \eqref{inequalityfornonlinear} that
\[
H'(t)\ge\frac{C_1^{1-p}}{2(1+t)^{(n-1)(p-1)/2}}H^p(t)
\quad\mbox{for}\ t\ge0,
\]
from which with the initial data $H(0)=(m(0)C_{0,g}/2)\e>0$
we can easily get the upper bound of lifespan estimate in Theorem \ref{blowup}.
\hfill $\Box$
\section{Proof of Theorem \ref{blowup1}}
\par\quad
The proof of Theorem \ref{blowup1} can be proceeded along almost the same way
as that of Theorem \ref{blowup}.
The only essential difference is that we replace a multiplier
$m$ defined in (\ref{test1}) by
\begin{equation}
\label{test1in}
\begin{aligned}
m_1(t):=(1+t)^{\mu},
\end{aligned}
\end{equation}
which was first introduced in Lai, Takamura and Wakasa \cite{LTW}
and has a property
\[
\frac{m'_1(t)}{m_1(t)}=\frac{\mu}{1+t}.
\]
Hence the differences in this section from the previous one
should appear only in points where
the boundedness of $m$ in (\ref{bound}) is employed.
They are (\ref{F1inequality}) and (\ref{inequalityfornonlinear}).
\par
Keeping this fact in our mind, we immediately obtain
\begin{lem}
\label{F1in}
Under the same assumption of Theorem \ref{blowup1}, it holds that
\begin{equation}
\label{F1inpostive}
F_1(t)\ge \frac{C_{f, 0}\e}{2m_1(t)}\ge0
\quad\mbox{for}\ t\ge0,
\end{equation}
where $F_1$ is defined in (\ref{F1def}).
\end{lem}
{\it Proof.} The proof is parallel to that of Lemma \ref{F1}.
Due to the unboundedness of $m_1$,
instead of (\ref{F1inequality}),
we have
\[
F_1'(t)+2F_1(t)\ge\frac{C_{f,g}\e}{m_1(t)}+\frac{1}{m_1(t)}\int_0^t\mu(1+s)^{\mu-1}F_1(s)ds
\]
by substituting $m$ with $m_1$ simply.
Integrating this inequality over $[0,t]$ with a multiplication $e^{2t}$, we get
\[
\begin{array}{ll}
e^{2t}F_1(t)\ge
&\d F_1(0)+C_{f,g}\e\int_0^t\frac{e^{2s}}{m_1(s)}ds\\
&\d+\int_0^t\frac{e^{2s}}{m_1(s)}ds\int_0^s\mu(1+r)^{\mu-1}F_1(r)dr.
\end{array}
\]
Therefore the comparison argument again yields that
\[
F_1(t)>\frac{C_{f,g}\e}{2m_1(t)}(1-e^{-2t})+e^{-2t}F_1(0)
\ge\frac{C_{f,0}\e}{2m_1(t)}\geq 0\quad\mbox{for}\ t\ge0
\]
as desired.
\hfill $\Box$
\vskip10pt
\par
In this way, we get the positivity of $F_1$ also for the case of $\beta=1$.
Due to this fact, we can proceed the proof of Theorem \ref{blowup1}
by simple replacement of $m$ by $m_1$
in the one of Theorem \ref{blowup}
till making use of the boundedness of $m$
once more at (\ref{inequalityfornonlinear}).
Hence, setting
\[
H_1(t):=
\frac{1}{2}\int_0^tm_1(s)ds\int_{\R^n}|u_t(x,s)|^p\psi(x,s)dx
+\frac{m_1(0)\e}{2}\int_{\R^n}g(x)\phi_1(x)dx,
\]
we have
\begin{equation}
\label{inequalityformpsiuin}
m_1(t)\int_{\R^n}u_t(x,t)\psi(x,t)dx\geq H_1(t) \quad\mbox{for}\ t\ge0.
\end{equation}
This is almost the same as (\ref{inequalityformpsiu}).
On the other hand, H\"{o}lder inequality and \eqref{inequalityforpsi}
as well as the concrete expression of $m_1$ yield that
\begin{equation}
\label{inequalityfornonlinearin}
2H_1'(t)\geq C_1^{1-p} (t+1)^{-(n+2\mu-1)(p-1)/2}
\left(m_1(t)\int_{\R^n}u_t(x,t)\psi(x,t)dx\right)^p.
\end{equation}
We then conclude from \eqref{inequalityformpsiuin} and \eqref{inequalityfornonlinearin} that
\[
H_1'(t)\ge\frac{C_1^{1-p}}{2(1+t)^{(n+2\mu-1)(p-1)/2}}H_1^p(t)
\quad\mbox{for}\ t\ge0,
\]
from which with the initial data $H_1(0)=(C_{0,g}/2)\e>0$
we can easily get the upper bound of lifespan estimate in Theorem \ref{blowup1}.
\hfill $\Box$

\begin{rem}
In the scale invariant damping case, we get an upper bound of the lifespan estimate depending on $\mu$,
since we have use the multiplier $m_1(t)=(1+t)^{\mu}$,
which is not bounded from above again, comparing to $m(t)$
used in the scattering case.
\end{rem}

\section*{Acknowledgment}
\par\quad
The first author is partially supported by NSFC(11501273), Zhejiang Province
Science Foundation(LY18A010008), Chinese Postdoctoral Science Foundation(2017M620128), NSFC(11771359,
11771194), high level talent project of Lishui City (2016RC25),
the Scientific Research Foundation of the First-Class Discipline of Zhejiang Province
(B)(201601). The second author is partially supported by the Grant-in-Aid for Scientific Research(C)
(No.15K04964),
Japan Society for the Promotion of Science,
and Special Research Expenses in FY2017, General Topics(No.B21),
Future University Hakodate.


\bibliographystyle{plain}

\begin{thebibliography}{20}


\bibitem{Agemi}
{R.Agemi}, {\it Blow-up of solutions to nonlinear wave equations in two space dimensions}, Manuscripta Math., {\bf 73} (1991), 153-162.

\bibitem{Glassey}
{R.T.Glassey}, MathReview to
\lq\lq
Global behavior of solutions to nonlinear wave equations
in three space dimensions"
of Sideris,
Comm. Partial Differential Equations (1983).

\bibitem{Hidano1}
{K.Hidano and K.Tsutaya}, {\it Global existence and asymptotic behavior of solutions for nonlinear wave
equations}, Indiana Univ. Math. J., {\bf 44} (1995), 1273-1305.

\bibitem{Hidano2}
{K.Hidano, C.Wang and K.Yokoyama}, {\it The Glassey conjecture with radially symmetric data}, J. Math.
Pures Appl., {\bf 98(9)} (2012), 518-541.


\bibitem{John1}
{F.John}, {\it Blow-up for quasilinear wave equations in three space dimensions}, Comm. Pure Appl. Math., {\bf 34} (1981), 29-51.
\bibitem{John2}
{F.John}, {\it Non-existence of global solutions of $\Box u=\frac{\partial}{\partial_t}F(u_t)$ in two and three space dimensions}, Rend. Circ. Mat. Palermo (2) Suppl., {\bf 8} (1985), 229-249.


\bibitem{LT}{N.-A.Lai and H.Takamura},
{\it Blow-up for semilinear damped wave equations
with sub-Strauss exponent
in the scattering case}, arXiv:1707.09583.


\bibitem{LTW}{N.-A.Lai, H.Takamura and K.Wakasa},
{\it Blow-up for semilinear wave equations
with the scale invariant damping and super-Fujita exponent},
J. Differential Equations, {\bf 263(9)} (2017), 5377-5394.


\bibitem{Masuda}
{K.Masuda}, {\it Blow-up solutions for quasi-linear wave equations in two space dimensions}, Lect. Notes Num. Appl. Anal., {\bf 6} (1983), 87-91.

\bibitem{Rammaha}{M.A.Rammaha},
{\it Finite-time blow-up for nonlinear wave equations in high dimensions},
Comm. Partial Differential Equations {\bf12(6)} (1987), 677-700.

\bibitem{Schaeffer}
{J.Schaeffer}, {\it Finite-time blow up for $u_{tt}-\Delta u=H(u_r, u_t)$ in two space dimensions}, Comm. Partial Differential Equations, {\bf 11 (5)} (1986), 513-543.


\bibitem{Sideris}
{T.C.Sideris}, {\it Global behavior of solutions to nonlinear wave equations in three space dimensions}, Comm. Partial Differential Equations, {\bf 8 (12)} (1983),
1291-1323.

\bibitem{Tzvetkov}
{N.Tzvetkov}, {\it Existence of global solutions to nonlinear massless Dirac system and wave equation
with small data}, Tsukuba J. Math., {\bf 22} (1998), 193-211.

\bibitem{Wirth1}{J.Wirth},
{\it Solution representations for a wave equation with weak dissipation},
Math. Methods Appl. Sci., {\bf 27} (2004),  101-124.

\bibitem{Wirth2}{J.Wirth},
{\it Wave equations with time-dependent dissipation. I. Non-effective dissipation},
J. Differential Equations, {\bf 222} (2006), 487-514.

\bibitem{Wirth3}{J.Wirth},
{\it Wave equations with time-dependent dissipation. II. Effective dissipation},
J. Differential Equations, {\bf 232} (2007),  74-103.

\bibitem{Yordanov}
{B.Yordanov and Q.S.Zhang}, {\it Finite time blow up for critical wave equations in high dimensions}, J. Funct. Anal., {\bf 231} (2006), 361-374.
\bibitem{Zhou1}
{Y.Zhou}, {\it Blow-up of solutions to the Cauchy problem for nonlinear wave equations}, Chin. Ann. Math., {\bf 22B (3)} (2001), 275-280.
\bibitem{Zhou2}
{Y.Zhou and W.Han}, {\it Blow-up of solutions to semilinear wave equations with variable coefficients and boundary}, J. Math. Anal. Appl., {\bf 374} (2011), 585-601.


\end{thebibliography}

\end{document}